\documentclass[12pt]{amsart}
\usepackage{graphicx}

\usepackage{amsthm}
\usepackage{amsmath}
\usepackage{amssymb}
\usepackage{url}
\usepackage[bookmarks=false]{hyperref}

\theoremstyle{plain}
\newtheorem{theorem}{Theorem}
\newtheorem{proposition}[theorem]{Proposition}

\theoremstyle{definition}

\newtheorem*{remark}{Remark}
\newtheorem{corollary}[theorem]{Corollary}
\newtheorem{lemma}[theorem]{Lemma}
\newtheorem*{example}{Example}
\newtheorem{question}{Question}
\newtheorem{conjecture}{Conjecture}

\newcommand{\Z}{\mathbb{Z}}
\newcommand{\G}{\mathcal{G}}

\title{
	{Equations of the Cayley Surface}
}
\author{Matty van Son}
\date{August 2021}

\thanks
{
	The author gratefully acknowledges partial support of the Austrian
	Science Fund (FWF) through project P 29981}


\begin{document}
	\begin{abstract}
		In this note we study the integer solutions of Cayley's cubic equation.
		We find infinite families of solutions built from recurrence relations.
		We use these solutions to solve certain general Pell equations.
		We also show the similarities and differences to Markov numbers.
		In particular we introduce new formulae for the solutions to Cayley's cubic equation in analogy with Markov numbers and discuss their distinctions.
	\end{abstract}
	
	\maketitle
	
	\section{Introduction}
	Cayley's equation $C_s(x,y,z)=0$ is defined by
	\[
	C_s(x,y,z)=s(x^2+y^2+z^2)-s^3-2xyz,
	\]
	and has been studied before in~\cite{cayley1869,hunt2000,heathbrown2003}.
	In this paper we study the positive integer solutions to Cayley's cubic equation.
	
	We have a new proof of the known fact that, when $s=1$, all solutions are given by Chebyshev polynomials.
	We find infinite families of solutions for each $s>1$ related to Chebyshev polynomials.
	The solutions are also shown to solve certain general Pell equations.
	Finally we compare the structure of solutions to the structure of Markov numbers, positive integer solutions to the equation $M(x,y,z)=0$, where
	\[
	M(x,y,z)=x^2+y^2+z^2-3xyz.
	\]
	
	\subsection{Background}
	This paper is motivated by the study of Markov numbers.
	It is well known that all solutions to the \textit{Markov Diophantine equation} $M(x,y,z)=0$ may be generated from any single solution $(a,b,c)$ by permutation and applying the transformation $(a,b,c)\mapsto (a,b,3ab-c)$.
	The books by T.~Cusick and M.~Flahive~\cite{cusick1989} and M.~Aigner~\cite{aigner2013} are good introductions to this area.
	
	Together with O.~Karpenkov~\cite{karpenkov2020} we have formulated a generalisation of Markov numbers.
	Our generalised Markov numbers are defined by recurrence relations of related continued fractions.
	We have developed an analogue to the following classical conjecture on Markov numbers.
	\begin{conjecture}[G.~F.~Frobenius, $1913$]
		Every solution to Markov's equation $(a,b,c)$ is uniquely defined by its largest element.
	\end{conjecture}
	The analogue of this \textit{uniqueness conjecture} fails for the generalised Markov numbers.
	It is known (see~\cite{mcginn2015}) that an analogue of the uniqueness conjecture also fails for the solutions to $C_1(x,y,z)=0$, which are also built by recurrence relations (see~\cite{denef1979,demeyer2007}).
	
	Our main aim is to study the solutions to $C_s(x,y,z)=0$ and compare the structure of the solutions with generalised Markov numbers.

	\subsection{Layout}

	\noindent In Section~\ref{section: cayley equations} we find all solutions to $C_s(x,y,z)=0$ for the $s=1$ case, and then find infinitely many solutions for the $s>2$ cases in Theorem~\ref{theorem: special recurrence}.
	We relate these to solutions of certain general Pell equations in Theorem~\ref{theorem: general pell}.
	
	\noindent In Section~\ref{section: s>1 results} we compare solutions to the Cayley equation with Markov numbers.
	We look briefly at the case $s=2$, and give some examples for $s>2$.

	\section{Cayley equations} \label{section: cayley equations}
	First note that Cayley's equation $C_s(x,y,z)=0$ is symmetric and quadratic in $x$, $y$, and $z$.
	Note also for any positive integer $p$ that we have infinitely many solutions 
	\[
	C_s(s,p,p)=s(s^2+2p^2)-s^3-2sp^2=0.
	\]
	
	\vspace{2mm}
	\noindent In Subsection~\ref{subsection: cheby} we find all solutions to $C_1(x,y,z)=0$ in Proposition~\ref{proposition: cheby soln full}, and use this to prove a well known result relating Pell equations to Chebyshev polynomials in Proposition~\ref{proposition: cheby pell}.
	
	\noindent In Subsection~\ref{subsection: general case} we study the more general $s>1$ cases, in particular we find infinite families of solutions in Theorem~\ref{theorem: special recurrence}.
	We use these solutions to solve certain Pell equations in Theorem~\ref{theorem: general pell}.
	
	\subsection{Solutions to $C_1(x,y,z)=0$.} \label{subsection: cheby}
	Let $s=1$.
	Then the equation $C_1(x,y,z)=0$ is solved by factorising 
	$ \displaystyle x=yz\pm\sqrt{(y^2-1)(z^2-1)}$.
	We have the following proposition.
	\begin{proposition} \label{proposition: conjugate}
		Suppose we have a solution $(a,b,c)$ to $C_1(x,y,z)=0$ with $b\geq\max\{a,c\}$ and $1<\min\{a,c\}$.
		Then
		\[
		\begin{aligned}
			b&=ac+\sqrt{(a^2-1)(c^2-1)},\\
			a&=bc-\sqrt{(b^2-1)(c^2-1)},\\
			c&=ab-\sqrt{(a^2-1)(b^2-1)}.
		\end{aligned}
		\]
	\end{proposition}
	
	\begin{proof}
		The second and third equalities hold simply since $a\leq b$ and $c\leq b$.
		Assume $b=ac-\sqrt{(a^2-1)(c^2-1)}$.
		We prove the first equality by contradiction. 
		As $c \leq b$ we have
		\[
		\begin{aligned}
			c&\leq b=ac-\sqrt{(a^2-1)(c^2-1)}\\
			\sqrt{(a^2-1)(c^2-1)}&\leq c(a-1)\\
			a^2c^2-a^2-c^2+1&\leq a^2c^2+c^2-2ac^2\\
			2c^2(a-1)&\leq a^2-1\\
			2c^2&\leq a+1.
		\end{aligned}
		\]
		The last inequality follows from $a>1$, all others follow from simple manipulations.
		Following the same argument for $a\leq b$ we have that $2a^2\leq c+1$.
		This is a contradiction, and so $b=ac+\sqrt{(a^2-1)(c^2-1)}$.
	\end{proof}
	
	Here and below we assume that any solution $C_1(a,b,c)=0$ has $b\geq \max\{a,c\}$ and $1<\min\{a,c\}$ unless otherwise stated.
	
	\begin{corollary}
		From a solution $C_1(a,b,c)=0$ we may generate three possibly distinct solutions, namely
		\[
		(\overline{a},b,c),\qquad (a,\overline{b},c),\qquad (a,b,\overline{c}),
		\]
		where $\overline{x}$ denotes the conjugate of $x$.
	\end{corollary}
	
	Since $\overline{b}\leq\max\{a,c\}$ the triple $(a,\overline{b},c)$ has a smaller maximal element than $(a,b,c)$.
	If we continue reducing the maximal element of triples in this way, eventually we will reach some triple $(x,y,z)$ with maximal element $x$ that cannot be reduced further.
	In this case $x=\overline{x}$ and so 
	\[
	(y^2-1)(z^2-1)=0.
	\]
	Hence $(x,y,z)=(x,x,1)$.
	We call a solution of this form a \textit{singular solution}.
	\begin{proposition} \label{proposition: quadratic reduction}
		Every solution to $C_1(x,y,z)=0$ may be reduced by repeated application of the quadratic formula to a singular solution.
	\end{proposition}
	
	\subsubsection{Lucas sequences}
	Let us define polynomial sequences that may be used to describe all solutions to $C_1(x,y,z)=0$.
	We later use these to polynomials in Section~\ref{section: s>1 results} to find infinitely many solutions to all $C_s(x,y,z)=0$.
	
	\textit{Lucas sequences of the first kind $U_n(P,Q)$} and \textit{of the second kind $V_n(P,Q)$} are the sequences defined by
	$U_0(P,Q)=0$, $U_1(P,Q)=1$, $V_0(P,Q)=2$, $V_1(P,Q)=P$, 
	and the recurrence relations
	\[
	\begin{aligned}
		U_{n+1}(P,Q)&=PU_n(P,Q)-QU_{n-1}(P,Q),\\
		V_{n+1}(P,Q)&=PV_n(P,Q)-QV_{n-1}(P,Q).
	\end{aligned}
	\]
	For example, $\big(U_n(1,{-}1)\big)$ are the Fibonacci numbers, and $\big(U_n(3,2)\big)$ are Mersenne numbers (numbers of the form $2^n{-}1$).

	\subsubsection{Cayley's cubic equation and Lucas sequences}
	Consider the polynomial sequences $T_n(x)=U_{n+1}(2x,1)$ and $\tilde{T}_n(x)=V_n(2x,1)$.
	These are called \textit{Chebyshev polynomials of the first kind} (respectively \textit{second kind}).
	The $T_n(x)$ are instrumental in generating solutions to $C_1(w,y,z)=0$.
	Later we use the $U_n(x)$ to solve a Pell equation related to $C_1(w,y,z)=0$.
	
	Note that any solution $C_1(1,p,p)=0$ for some positive integer may be recast as a solution in terms of Chebyshev polynomials as $C_1\big(T_0(p),T_1(p),T_1(p)\big)=0$.
	\begin{proposition} \label{proposition: cheby soln}
		Assume that a solution to $C_1(x,y,z)=0$ may be written in terms of some $p$ by
		\[
		C_1\big(x,T_n(p),T_m(p)\big)=0,
		\]
		for some positive integers $n$ and $m$.
		Then
		\[
		x\in\{T_{n+m}(p),T_{|n-m|}(p)\}.
		\]
	\end{proposition}

	To prove Proposition~\ref{proposition: cheby soln} we need the following two lemmas.
	The first is a classical result on Chebyshev polynomials which we state without proof.
	\begin{lemma} 
		We have that $\displaystyle 2T_n(p)T_m(p)=T_{n+m}(p)+T_{|n-m|}(p)$.
	\end{lemma}
	
	The following lemma gives a formula for $T_{n+m}(p)-T_{|n-m|}(p)$.
	\begin{lemma} \label{lemma: cheby square root}
		We have that 
		\[
		4\big(T_n(p)^2-1\big)\big(T_m(p)^2-1\big)=\big(T_{n+m}(p)-T_{|n-m|}(p)\big)^2.
		\]
	\end{lemma}
	
	\begin{proof}
		We shorten $T_n(p)$ to $T_n$ in this proof.
		From the previous lemma we have that $T_n^2=T_{2n}/2+1$.
		Assuming, without loss of generality, that $m\leq n$, then
		\[
		\begin{aligned}
			8\big(T_n^2-1\big)\big(T_m^2-1\big)&=2\big(T_{2n}-1\big)\big(T_{2m}^2-1\big)\\
			&=2\big(T_{2n}T_{2m}-T_{2n}-T_{2m}+1\big)\\
			&=T_{2n+2m}+T_{2n-2m}-2T_{2n}-2T_{2m}+2.
		\end{aligned}
		\]
		We also have that
		\[
		\begin{aligned}
			2\big(T_{n+m}-T_{|n-m|}\big)^2&=2T_{n+m}^2-4T_{n+m}T_{n-m}+2T_{n-m}^2\\
			&=T_{2n+2m}+T_{2n-2m}+2-2T_{2n}-2T_{2m}.
		\end{aligned}
		\]
		The statement holds for all $n$ and $m$ and the proof is complete.
	\end{proof}
	
	\begin{proof}[Proof of Proposition~\ref{proposition: cheby soln}]
		We shorten $T_n(p)$ to $T_n$ in this proof.
		Factorising $C_1\big(x,T_n(p),T_m(p)\big)=0$ we have
		\[
		\begin{aligned}
			x&=T_nT_m\pm\sqrt{\big(T_n^2-1\big)\big(T_m^2-1\big)}\\
			&=\frac{T_{n+m}+T_{|n-m|}}{2}\pm\frac{\big(T_{n+m}-T_{|n-m|}\big)}{2}
			&\in\{T_{n+m},T_{|n-m|}\}.
		\end{aligned}
		\]
		The second equality follows from the previous lemmas.
		This completes the proof.
	\end{proof}
	
	From the base solutions $C_1\big(T_0(p),T_1(p),T_1(p)\big)=0$ we can generate, through the quadratic formula, solutions
	\[
	C_1\big(T_n(p),T_{n+m}(p),T_m(p)\big)=0
	\]
	for any positive integers $n$ and $m$ with $\gcd(n,m)=1$.
	This is a consequence of the study of Euclid trees, to which one can find references with a similar application in~\cite[p.~211-212]{aigner2013} and~\cite[p.~15-22]{mcginn2015}.
	
	Let us recall another classical result on the composition of Chebyshev polynomials.
	\begin{lemma}
		For any positive integers $a$ and $b$ we have
		\[
		T_{ab}(p)=T_a\big(T_b(p)\big).
		\]
	\end{lemma}
	As a corollary to this result, if $\gcd(n,m)>1$ and
	\[
	n=\gcd(n,m)t, \quad
	m=\gcd(n,m)s, \quad
	q=T_{\gcd(n,m)}(p),
	\]
	for some positive integers $t$ and $s$, with $\gcd(s,t)=1$, then 
	\[
	C_1\big(T_n(p),T_{n+m}(p),T_m(p)\big)=C_1\big(T_t(q),T_{t+s}(q),T_s(q)\big)=0.
	\]
	Hence we have the following known result.
	\begin{proposition} \label{proposition: cheby soln full}
		All positive integer solutions to $C_1(x,y,z)=0$ are given by
		\[
		\Big(T_n(p),T_{n+m}(p),T_m(p)\Big)
		\]
		where $p$ is any positive integer and $n$ and $m$ are non negative integers.
	\end{proposition}
	
	\begin{proof}
		By Proposition~\ref{proposition: quadratic reduction} every solution $(a,b,c)$ reduces by the quadratic formula to $(1,x,x)=\big(T_0(x),T_1(x),T_1(x)\big)$ for some positive integer $x$.
		Reversing this reduction using Proposition~\ref{proposition: cheby soln} we see that there are non negative integers $n$ and $m$ with $(a,b,c)=\big(T_n(x),T_{n+m}(x),T_m(x)\big)$.
	\end{proof}
	
	\begin{remark}
		One can find values of $x$ for which
		\[
		x=yz\pm\sqrt{(y^2-1)(z^2-1)}\in\Z
		\]
		by assuming (without loss of generality) that $y\leq z$.
		If $z^2-1=(y^2-1)a^2$ for some positive integer $a$ then 
		\[
		x=yz\pm (y^2-1)a\in\Z.
		\]
		This leads to solving the Pell equation $z^2-da^2=1$,
		where $d=y^2-1$ and $z$ and $a$ are to be found as functions of $y$.
		The above theorem may be restated in the following way: all positive integer solutions to
		\[
		z^2-(y^2-1)a^2=1
		\]
		are given by $z=T_n(y)$ and $a=U_{n-1}(y)$
		for some positive integer $n$.
		This is a known result and serves as an alternate definition for Chebyshev polynomials, see~\cite{denef1979,demeyer2007}.
	\end{remark}
	
	\begin{example}
		Let $y=2$.
		Then the solutions $(z,a)$ to $z^2-3a^2=1$ begin
		\[
		\Big(\big(T_n(3),U_{n-1}(3)\big)\Big)=\big((2,1),(7,4),(26,15),\ldots,(1351,780),\ldots\big).
		\]
		The fraction $1351/780$ has been used as an approximation to $\sqrt{3}$ since antiquity.
	\end{example}
	
	The next proposition solves a similar Pell equation, and is a rewording of Lemma~\ref{lemma: cheby square root}.
	\begin{proposition} \label{proposition: cheby pell}
		Let $d=y^2-1$.
		Let $y=T_n(p)$ for some positive integers $n$ and $p$.
		Then the Pell equation 
		\[
		a^2-dz^2=-d
		\]
		is also solved for $z$ and $a$, and some positive integer $m$, by
		\[
		z=T_m(p), \quad a=\frac{T_{n+m}(p)-T_{|n-m|}(p)}{2}.
		\]
	\end{proposition}
	
	\begin{proof}
		From Lemma~\ref{lemma: cheby square root} we have that 
		\[
		a^2=\big(T_n(p)^2-1)(T_m(p)^2-1\big)=d(T_m(p)^2-1\big).
		\]
		Hence
		\[
		a^2-dz^2=d\big(T_m(p)^2-1\big)-dT_m(p)^2=-d
		\]
	\end{proof}
	
	\begin{example}
		Let $n=3$ and $p=4$.
		Then $d=T_3(4)^2-1=960$.
		Hence the solutions to $a^2-960z^2=-960$ begin
		\[
		\left( \frac{T_{3+m}(4)-T_{|3-m|}(4)}{2},T_m(4) \right)=\big((120,4),(960,31),(7560,244),\ldots\big).
		\]
	\end{example}

	\subsection{Solutions to $C_s(x,y,z)=0$} \label{subsection: general case}

	In this subsection we find infinite families of positive integer solutions to the equation $C_s(x,y,z)=0$.
	We show that these solutions are generated by polynomial sequences, similar to how Chebyshev polynomials generate solutions when $s=1$.
	There are still solutions not constructed in this way, we give an example of this later in Subsection~\ref{subsection: Cs examples}.
	
\begin{figure}
		\centering\includegraphics{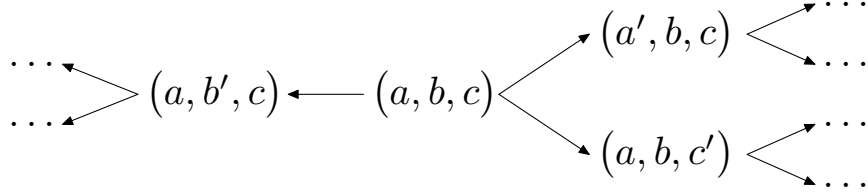}
		\caption{The graph $\G(a,b,c)$.} \label{figure: graph of cayley}
	\end{figure}

	Let $(a,b,c)$ be a solution to $C_s(x,y,z)=0$, with no assumptions on the size of elements.
	Assume that $s$ divides one of the elements of this solution, say $s$ divides $b$.
	Let $\overline{a}$ and $\overline{c}$ denote the conjugates of $a$ and $c$ respectively.
	Then $\overline{a}$ and $\overline{c}$ are both positive integers (The argument is similar to the proof of Proposition~\ref{proposition: conjugate}).
	As such, there is a graph $\G(a,b,c)$ (see Figure~\ref{figure: graph of cayley}) whose vertices are solutions to $C_s(x,y,z)=0$,
	where two vertices $v=(v_1,v_2,v_3)$ and $w$ are connected by an edge $(v,w)$ if $w$ is one of the following
	\[
	(\overline{v_1},v_2,v_3),\quad (v_1,\overline{v_2},v_3), \quad (v_1,v_2,\overline{v_3}).
	\]

	\begin{example}
		Let $s=3$ and $b=6$.
		Then $C_3(21, 4053, 291)=0$.
		We have that
		\[
		\begin{aligned}
			\frac{21\cdot 4053+\sqrt{(21^2-3^2)(4053^2-3^2)}}{3}&=56451,\\
			\frac{291\cdot 4053+\sqrt{(291^2-3^2)(4053^2-3^2)}}{3}&=786261,\\
			\frac{21\cdot291-\sqrt{(21^2-3^2)(291^2-3^2)}}{3}&=21,\\
		\end{aligned}
		\]
		and indeed $\displaystyle C_3(21, 56451, 4053)=0$ and $\displaystyle C_3(4053, 786261, 291)=0$.
	\end{example}
	
	\begin{remark}
		There is a degenerate case of $\overline{a}=a$ (or $\overline{c}=c$), where the triple $(a,b,c)$ is in fact the triple $(s,p,p)$ for some positive integer $p>s$.
		In this case we can only generate different solutions from the quadratic formula if $s$ divides $2p^2$.
		This is clear since
		\[
		\overline{s}=\frac{p^2+\sqrt{(p^2-s^2)^2}}{s}=\frac{2p^2-s^2}{s}.
		\] 
		Moreover, from the triple $(\overline{s},p,p)$ we can only generate new solutions if 
		\[
		\overline{p}=\frac{\overline{s}p+\sqrt{(\overline{s}^2-s^2)(p^2-s^2)}}{s}
		\]
		is also an integer.
		If $s$ does not divide $2p^2$ then $\G(s,p,p)$ contains a single vertex.
		If $s$ divides $2p^2$ but $\overline{p}\notin \Z$ then $\G(s,p,p)$ contains the two vertices $(s,p,p)$ and $(\overline{s},p,p)$, connected by an edge.
		In all other cases we generate an infinite graph, as detailed below.
	\end{remark}

	Let $v=(a,b,c)$ be a vertex in $\G(a,b,c)$.
	If $\overline{a}\neq a$ and $\overline{c}\neq c$ then we have two new vertices $v_1=(\overline{a},b,c)$ and $w_1=(a,b,\overline{c})$ in the graph $\G(a,b,c)$, both connected to $v$ by the edges $(v,v_1)$ and $(v,w_1)$.
	Reapplying this argument (based on the fact that $s$ divides $b$) we can generate infinitely many solutions and build a binary graph.
	\begin{example}
		When $s=1$ or $s=2$ then $\G(s,p,p)$ is an infinite graph for any $p> 1$.
	\end{example}

	The following proposition details the relation between consecutively generated solutions.
	\begin{proposition} \label{proposition: rec rel general}
		Consider a solution $(a,b,c)$ to $C_s(x,y,z)=0$ with $b\geq\max\{a,c\}$.
		Assume without loss of generality that $a\leq c$, and assume that $\overline{a}$ and $\overline{c}$ are integers.
		Then 
		\[
		\overline{a}=\frac{2b}{s}c-a, \quad
		\overline{c}=\frac{2b}{s}\overline{a}-c.
		\]
	\end{proposition}

	\begin{proof}
		We are assuming that neither $\overline{a}$ nor $c$ are zero. 
		Note that
		\[
		\overline{c}=\frac{\overline{a}b+\sqrt{(\overline{a}^2-s^2)(b^2-s^2)}}{s},\qquad c=\frac{\overline{a}b-\sqrt{(\overline{a}^2-s^2)(b^2-s^2)}}{s}
		\]
		and so
		\[
		\frac{\overline{c}+c}{\overline{a}}=\frac{2b}{s}.
		\]
		Similarly, using the formula for $\overline{a}$ and $a$ we have
		\[
		\frac{\overline{a}+a}{c}=\frac{2b}{s}.
		\]
	\end{proof}
	
	\subsubsection{Lucas sequences when $s>1$}
	
	Let us define a sequence $R_{n,s}(x)$ to help solve $C_s(w,y,z)=0$.
	We also define an alternate sequence $R_{n,s}^*(x)$ that we use later to solve certain Pell equations.
	Although the sequences are dependant on $s$ we will write $R_n(x)$ and $R_n^*(x)$ where $s$ is clear.
	We define $R_n(x)$ and $R_n^*(x)$ in terms of Lucas sequences
	\[
	R_n(x)=\frac{s}{2}V_n\left(\frac{2x}{s},1\right),\qquad R_n^*(x)=U_{n+1}\left(\frac{2x}{s},1\right).
	\]
	In particular, $R_0(x)=s$, $R_1(x)=x$, $R_0^*(x)=1$, and $R_1^*(x)=2x/s$.
	We relate $R_n(x)$ to $T_n(x)$.

	\begin{proposition}
		We have that
		\[
		R_n(x)=sT_n\left(\frac{x}{s}\right).
		\]
	\end{proposition}
	
	\begin{proof}
		As $R_0(x)=s=1s=sT_0(x/s)$ and
		\[
		R_1(x)=x=\frac{sx}{s}=sT_1(x/s)
		\]
		we use with induction.
		Assume that $R_n(x)=sT_n(x/s)$ for $n=0,\ldots,k$, where $k$ is some positive integer.
		Then
		\[
		\begin{aligned}
			R_{k+1}(x)&=\frac{2x}{s}R_k(x)-R_{k-1}(x)\\
			&=2xT_k\left(\frac{x}{s}\right)-sT_{k-1}\left(\frac{x}{s}\right)\\
			&=s\left(\frac{2x}{s}T_k\left(\frac{x}{s}\right)-T_{k-1}\left(\frac{x}{s}\right)\right)=sT_{k+1}\left(\frac{x}{s}\right).
		\end{aligned}
		\]
		The second equality follows from the induction hypothesis.
		This completes the induction and the proof.
	\end{proof}
	
	\begin{remark}
		Note that $C_s(b,\overline{a},c)=0$ and $C_s(b,\overline{c},\overline{a})=0$.
		Then 
		\[
		C_s\big(b,R_i(b),R_{i-1}(b)\big)=0
		\]
		for all $i>1$.
		All $R_n(b)$ are integers if and only if $s$ divides $2b$.
	\end{remark}

	\subsubsection{Solutions to $C_s(x,y,z)=0$}
	We have the following theorem.
	\begin{theorem} \label{theorem: special recurrence}
		Let $s$ and $b$ be positive integers where $s$ divides $2b$.
		Then
		\[
		C_s\big(R_n(b),R_{n+m}(b),R_m(b)\big)=0
		\]
		for any non negative integers $n$ and $m$. 
	\end{theorem}
	
	\begin{remark}
		Note that
		\[
		R_{n+1}(b)=\frac{bR_n(b)+\sqrt{\big(R_n(b)^2-s^2\big)\big(b^2-s^2\big)}}{s},
		\]
		and, by definition,
		\[
		C_s\big(b,R_n(b),R_{n-1}(b)\big)=0
		\]
		for all $n\geq 1$.
	\end{remark}

	\begin{proof}[Proof of Theorem~\ref{theorem: special recurrence}]
		Let $x=b/s$.
		Recall that
		\[
		T_n(x)^2+T_{n+m}(x)^2+T_m(x)^2=1+2T_n(x)T_{n+m}(x)T_m(x).
		\]
		Then
		\[
		\begin{aligned}
			&C_s\big(R_n(b),R_{n+m}(b),R_m(b)\big)\\
			=&C_s\big(sT_n(x),sT_{n+m}(x),sT_m(x)\big)\\
			=&s^3\big(T_n(x)^2+T_{n+m}(x)^2+T_m(x)^2\big)-s^3-2s^3T_n(x)T_{n+m}(x)T_m(x)\\
			=&0.
		\end{aligned}
		\]
		This completes the proof.
	\end{proof}
	
	Although we can not find all solutions in this way, as we can when $s=1$, we can still analyse solutions if they reduce to a singular $(s,b,b)$.
	
	\begin{corollary}
		Consider a solution $(v_1,v_2,v_3)$ to $C_s(x,y,z)=0$.
		Assume $(v_1,v_2,v_3)$ may be reduced by the quadratic formula to a solution $C_s(s,b,b)=0$ for some positive integer $b$.
		Assume also that $s$ divides $2b$.
		Then there exist non negative integers $n$ and $m$ such that
		\[
		\big\{v_1,v_2,v_3\big\}=\big\{R_n(b),R_{n+m}(b),R_m(b)\big\}.
		\]
	\end{corollary}
	
	\begin{remark}
		In particular, if $C_s(b,x,y)=0$ and $s$ divides $2b$, then for some positive integer $n$ we have that $x=R_n(b)$ and $y=R_{n\pm1}(b)$.
	\end{remark}
	
	Finally for this section, let us use these results to tackle certain Pell equations.
	
	\begin{theorem} \label{theorem: general pell}
		Let $d=y^2-4$.
		
		\begin{itemize}
			\item[(\emph{i})] 
			The Pell equation
			\[
			z^2-da^2=s^2
			\]
			is solved for $z$ and $a$ by
			\[
			z=R_n(y), \quad a=R_n^*(y).
			\]
			\item[(\emph{ii})] 
			Let $y=R_n(p)$ for some positive integers $n$ and $p$.
			Then the Pell equation 
			\[
			a^2-dz^2=-s^2d
			\]
			is solved for $z$ and $a$ by
			\[
			z=R_m(p), \quad a=R_{n+m}(p)-R_{|n-m|}(p).
			\]
		\end{itemize}

	\end{theorem}
	
	\begin{proof}
		We shorten $R_n(y)$ to $R_n$ in this proof.
		For $n\geq 2$ we prove the following two statements together:
		\begin{itemize}
			\item $\displaystyle R_n^2-d(R_{n-1}^{*})^2=s^2$;
			\item $\displaystyle R_nR_{n-1}-dR_{n-1}^*R_{n-2}^*=sy$.
		\end{itemize}
		
		For $n=1$ we have that $y^2-(y^2-s^2)=s^2$.
		When $n=2$ we have 
		\[
		(4y^4/s^2+s^2-4y^2)-(y^2-s^2)(4y^2/s^2)=s^2
		\]
		and $2y^3/s-sy-2y^3/s+2sy=sy$.
		Hence both statements hold for $n=2$.
		Assume both statements are true for all $n=2,\ldots,k$, some positive integer $k>2$.
		Then
		\[
		\begin{aligned}
			&R_kR_{k+1}-dR_k^*R_{k-1}^*\\
			=&\frac{2y}{s}\big(R_k^2-d(R_{k-1}^*)^2\big)-R_kR_{k-1}+dR_{k-1}^*R_{k-2}^*\\
			=&2sy-sy=sy.
		\end{aligned}
		\]
		The first equality substitutes in the recurrence relations for $R_{k+1}$ and $R_k^*$, while the second equality follows from the induction hypotheses.
		Next we have
		\[
		\begin{aligned}
			&R_{k+1}^2-(y^2-s^2)(R_{k+1}^*)^2\\
			=&\frac{4y^2}{s^2}\Big(R_k^2-d(R_{k-1}^*)^2\Big) +R_{k-1}^2-d(R_{k-2}^*)^2-\frac{4y}{s}\big(R_kR_{k-1} {-} dR_{k-1}^*R_{k-2}^*\big)\\
			=&4y^2+s^2-4y^2=s^2.
		\end{aligned}
		\]
		The first equality is again found by substituting the  recurrence relations for $R_{k+1}$ and $R_{k+1}^*$, while the second and third equalities follow from the induction hypotheses.
		
		Part (\emph{ii}) follows from the relation
		\[
		4\big(R_n(p)^2-s^2\big)\big(R_m(p)^2-s^2\big)=s^2\big(R_{n+m}(p)-R_{|n-m|}(p)\big)^2,
		\]
		which in turn follows from Lemma~\ref{lemma: cheby square root} with substitution
		\[
		R_n(p)=sT_n(p/s).
		\]
		This completes the proof.
	\end{proof}
	
	\section{General Markov numbers and experimental results for $C_s(x,y,z)=0$} \label{section: s>1 results}
	In Subsection~\ref{subsection: markov} we show that no solutions that we have found to $C_s(x,y,z)=0$ represent general Markov numbers .
	We make note of the sequences and relations for the case $s=2$ in Subsection~\ref{subsection: lucas} .
	Then in Subsection~\ref{subsection: Cs examples} we give some experimental results for $s>2$, and make note of some unresolved questions in these cases.
	
	\subsection{Markov recurrence relations} \label{subsection: markov}
	In this subsection we show that the recurrence relation for general Markov numbers is different to the recurrence relation for the polynomials $R_n(x)$.
	
	\subsubsection{Continuant definitions}
	We start with the necessary definitions for general Markov numbers, starting with continued fractions.
	
	Any real number $r$ may be written in the form
	\[
	r=a_0+\cfrac{1}{a_1+\cfrac{1}{\ddots}},
	\]
	where $a_0\in\Z$ and $a_1,a_2,\ldots$ are positive integers.
	The right hand side of this equation is called a \textit{continued fraction} and we denote it by the brackets $[a_0;a_1:a_2:\ldots]$.
	
	Let $r=[a_0;a_1,\ldots:a_n]$ for some positive integer $n$.
	Then $r=p/q$ for some positive integers with $\gcd(p,q)=1$.
	We call $p$ the \textit{continuant} of the sequence $(a_0,\ldots,a_n)$, denoted by $K(a_0,\ldots,a_n)$.
	Continuants may be alternatively defined by $K()=1$, $K(x_0)=x_0$, and the recursive formula
	\[
	\begin{aligned}
		K(x_0,\dots,x_m)&=x_0K(x_1,\ldots,x_m)+K(x_2,\ldots,x_m).
	\end{aligned}
	\]
	\begin{remark}
		Let $\alpha=(a_0,\ldots,a_n)$.
		We use the following notation for abbreviated continuants
		\[
		\begin{aligned}
			\breve{K}(a_0,a_1,\ldots,a_n)&=\breve{K}(\alpha)=K(a_0,a_1,\ldots,a_{n-1}),\\
			\breve{K}_2(a_0,a_1,\ldots,a_n)&=\breve{K}_2(\alpha)=K(a_1,\ldots,a_{n-1}).
		\end{aligned}
		\]
		Let $\beta=(b_0,\ldots,b_m)$.
		Then we write $\alpha\beta=(a_0,\ldots,a_n,b_0,\ldots,b_m)$ and $\alpha^k=\alpha\alpha\ldots$ 
		where $\alpha$ is repeated $k$ times.
		The following is particular version of a classical formula for splitting continuants (one can find a proof in~\cite{graham1990}).
		\begin{equation} \label{equation: cont split}
		\breve{K}(\alpha^2\beta)=K(\alpha)\breve{K}(\alpha\beta)+\breve{K}(\alpha)\breve{K}_2(\alpha\beta).
		\end{equation}
	\end{remark}
	
	In~\cite{karpenkov2020} we have found the following recurrence relation for continuants of even length sequences.
	All general Markov numbers are defined by sequences that must satisfy this recursion.
	\begin{proposition}
		Let $\alpha=(a_1,\ldots,a_{2m})$, $\lambda$, and $\rho$ be sequences of positive integers.
		Then
		\[
		\frac{\breve{K}(\alpha^2)}{\breve{K}(\alpha)}=\frac{\breve{K}(\lambda\alpha^2\rho)+\breve{K}(\lambda\rho)}{\breve{K}(\lambda\alpha\rho)}.
		\]
	\end{proposition}

	\begin{remark}
		In particular, for even length sequences $\alpha$ and $\beta$ we have
		\[
			\breve{K}(\alpha^{k+1}\beta)=\frac{\breve{K}(\alpha^2)}{\breve{K}(\alpha)}\breve{K}(\alpha^k\beta)-\breve{K}(\alpha^{k-1}\beta),
		\]
		The sequence built from recursion is $\displaystyle \big(\breve{K}(\beta),\breve{K}(\alpha\beta),\breve{K}(\alpha^2\beta),\ldots\big)$.
		General Markov numbers occur in triples
		\[
		\Big(\breve{K}(\alpha),\breve{K}(\alpha^n\beta),\breve{K}(\alpha^{n+1}\beta)\Big).
		\]
		When $\big(\breve{K}(\alpha),\breve{K}(\alpha\beta),\breve{K}(\beta)\big)$ is a vertex in the graph of regular Markov numbers, then
		\[
		\breve{K}(\alpha)^2+\breve{K}(\alpha^n\beta)^2+\breve{K}(\alpha^{n-1}\beta)^2=3\breve{K}(\alpha\beta)\breve{K}(\alpha^n\beta)\breve{K}(\alpha^{n-1}\beta).
		\]
		Similar formulae exist when $\beta$ is repeated.
	\end{remark}

	\subsubsection{Recurrence relation comparison}
	All infinite families of solutions to $C_s(x,y,z)=0$ that we have found are related to the sequence of polynomials $R_n(b)$,
	where $s$ divides $2b$.
	Recall from the Cayley equation that 
	\[
	C_s\big(R_1(b),R_n(b),R_{n+1}(b)\big)=0.
	\]
	To align with general Markov numbers we assume that
	\[
	\breve{K}(\alpha)=R_1(b), \quad \breve{K}(\alpha\beta)=R_2(b),\quad \frac{\breve{K}(\alpha^2)}{\breve{K}(\alpha)}=\frac{2b}{s}.
	\]
	(The following argument also works when $\beta$ is the repeated sequence.)
	From the recursive formulae for both the continuants and $R_n(b)$ we find that $R_0(b)=\breve{K}(\beta)$.
	However this implies that
	\[
	\frac{2\breve{K}(\alpha)}{\breve{K}(\beta)}=\frac{2b}{s}=\frac{\breve{K}(\alpha^2)}{\breve{K}(\alpha)}=K(\alpha)+\breve{K}_2(\alpha),
	\]
	and in particular that $2\breve{K}(\alpha)>K(\alpha)\breve{K}(\beta)$.
	This can only happen when $s=\breve{K}(\beta)=1$.
	In this case we have that
	\[
	T_n(b)=\breve{K}(\alpha^{n-1}\beta).
	\]
	Applying the product formula for Chebyshev polynomials we get that
	\[
	2\breve{K}(\alpha^n\beta)\breve{K}(\alpha\beta)=\breve{K}(\alpha^{n+1}\beta)+\breve{K}(\alpha^{n-1}\beta),
	\]
	which, once compared to the recurrence relation for continuants, gives
	\[
	K(\alpha)+\breve{K}_2(\alpha)=\frac{\breve{K}(\alpha^2)}{\breve{K}(\alpha)}=2\breve{K}(\alpha\beta).
	\]
	Since $\breve{K}(\beta)=1$, either $\beta=(x)$ or $\beta=(1,x)$ for some positive integer $x$.
	If $\beta=(1,x)$ then
	\[
	K(\alpha)+\breve{K}_2(\alpha)=2K(\alpha)+2\breve{K}(\alpha),
	\]
	which can not be true.
	If $\beta=(x)$ then 
	\[
	K(\alpha)+\breve{K}_2(\alpha)=2K(\alpha)
	\]
	which again can not be true.	
	Hence no sequence of general Markov numbers is given by the polynomial sequence $R_n(x)$.

	\subsection{Solutions to $C_2(x,y,z)=0$} \label{subsection: lucas}
	Note for any positive integers $t$ and $s$ that
	\[
	\frac{t+s}{2}\in\Z\quad\mbox{if and only if}\quad\frac{t-s}{2}\in\Z.
	\]
	Hence all solutions to $C_2(x,y,z)=0$ may be given
	\[
	C_2\big(R_{n,2}(p),R_{n+m,2}(p),R_{m,2}(p)\big)=0.
	\]
	It is simple to check that $R_{n,2}(p)=V_n(p,1)$ and $R_{n-1,2}^*(p)=U_n(p,1)$.
	\begin{remark}
		We write $L_n(p)=V_n(p,1)$ and $F_n(p)=U_n(p,1)$.
		The $L_n$ and $R_n$ are closely related to Lucas and Fibonacci polynomials respectively. 
		The difference is in the recurrence relation, which are built by addition rather than subtraction.
		
		The $L_n(x)$ have been called \textit{auxiliary Rudin-Shapiro polynomials}, see sequence A$213234$ on the OEIS~\cite{oeis}.
		The relation between $L_n(x)$ and $T_n(x)$ is 
		\[
		L_n\left(2x\right)=2T_n\left(x\right).
		\]
	\end{remark}
	
	\begin{corollary}
		Let $d=y^2-4$.
		\begin{itemize}
			
			\item[(\emph{i})] The Pell equation $\displaystyle z^2-da^2=4 $
			is solved for $z$ and $a$ by
			\[
			z=L_n(y), \quad a=F_n(y).
			\]
			
			\item[(\emph{ii})] Let $y=L_n(p)$ for some positive integers $n$ and $p$.
			Then the Pell equation $\displaystyle a^2-dz^2=-4d$
			is also solved for $z$ and $a$ by
			\[
			z=L_m(p), \quad a=L_{n+m}(p)-L_{|n-m|}(p),
			\]
			for some positive integer $m$.
		\end{itemize}

	\end{corollary}

	\subsection{Some examples for $C_s(x,y,z)=0$} \label{subsection: Cs examples}
	We have said that not all solutions to $C_s(x,y,z)=0$ are given by the formula
	\[
	C_s\big(R_n(p),R_{n+m}(p),R_m(p)\big)=0.
	\]
	The following examples show this.
	
	\begin{example}
		There exist isolated solutions $C_s(a,b,c)=0$ where 
		$\overline{a}$, $\overline{b}$, and $\overline{c}$ are all non integer.
		Let $s=24$ and $(a,b,c)=(26,74,51)$. 
		Then
		\[
		C_{24}(26,74,51)=0,\quad \overline{a}=\frac{577}{2}, \quad \overline{b}=\frac{73}{2}, \quad \overline{c}=\frac{328}{3}.
		\]
		
		Sometimes solutions appear in pairs, as in the following.
		\[
		\begin{aligned}
			C_{12}(13,20,15)&=0,\\
			C_{12}(37,20,15)&=0.
		\end{aligned}
		\]
		Neither of these solutions generate any others besides themselves.
	\end{example}
	We finish with an open question about the solutions $C_s(x,y,z)=0$.
	\begin{question}
		Can a formula be found to describe all solutions to $C_s(x,y,z)=0$?
		Can it be shown that such a formula does not exist?
	\end{question}

	\bibliographystyle{plain}
	\bibliography{biblio_ALL}

\end{document}